%% file: main.tex
\documentclass[letterpaper]{article} % DO NOT CHANGE THIS
\usepackage[preprint]{aaai2027}  % DO NOT CHANGE THIS
\usepackage[hyphens]{url}  % DO NOT CHANGE THIS
\usepackage{graphicx} % DO NOT CHANGE THIS
\usepackage{natbib}  % DO NOT CHANGE THIS AND DO NOT ADD ANY OPTIONS TO IT
\usepackage{caption} % DO NOT CHANGE THIS AND DO NOT ADD ANY OPTIONS TO IT
\usepackage[table]{xcolor}
\usepackage{algorithm}
\usepackage{amsmath}
\usepackage{mathrsfs}
\usepackage{mathtools}
\usepackage{amsfonts}
\usepackage{algorithmic}
\usepackage{enumitem}
\newcommand{\R}{\mathbb{R}}
\newcommand{\E}{\mathbb{E}}

\newcommand{\argmin}{\mathop{\arg\min}}
\newcommand{\lra}[1]{\left(#1\right)}
\newcommand{\lrb}[1]{\left[#1\right]}
\newcommand{\lrc}[1]{\left\{#1\right\}}

\usepackage{newfloat}
\usepackage{listings}
\DeclareCaptionStyle{ruled}{labelfont=normalfont,labelsep=colon,strut=off} % DO NOT CHANGE THIS
\floatstyle{ruled}
\newfloat{listing}{tb}{lst}{}
\floatname{listing}{Listing}

\usepackage{booktabs}

\title{BOCO: Bayesian Online Contextual Optimization\\for Decision-Focused Online Learning}
\author{
    Zhuojun Xie,
    Adam Abdin,
    Yiping Fang\corresponding,
}
\affiliations{
    \textsuperscript{\rm 1}Laboratoire G\'{e}nie Industriel, CentraleSup\'{e}lec, Universit\'{e} Paris-Saclay\\
    3 Rue Joliot Curie, 91190 Gif-sur-Yvette, France\\

    xie.zhuojun, adam.abdin, yiping.fang@centralesupelec.fr
}

\begin{document}
\maketitle

\begin{abstract}
\emph{Decision-focused learning} (DFL) trains predictive models to optimize downstream decisions rather than prediction accuracy alone.
While recent studies have extended this paradigm to online settings with streaming data, existing online DFL methods generally maintain a point estimate, while their gradient-based updates require either a differentiable optimization layer or a problem-specific surrogate loss. Consequently, they can be unstable under limited data and difficult to apply across heterogeneous optimization problems.
We introduce Bayesian Online Contextual Optimization (\texttt{BOCO}), a framework that maintains a decision-focused posterior over model parameters. \texttt{BOCO} aggregates the resulting predictions when prescribing decisions, thereby accounting for parameter uncertainty.
To track this posterior in evolving environments, we develop two particle-based inference algorithms: a sequential Monte Carlo sampler for nondifferentiable problems and a function-space Stein variational gradient descent algorithm for differentiable problems. Across both real-world tasks, \texttt{BOCO} reduces running mean regret and temporal regret variability relative to two frequentist online DFL baselines. At the full horizon, relative to the best frequentist baseline, BOCO reduces running mean regret and temporal regret variability by 6.9\% and 6.3\% on knapsack and by 49.3\% and 38.1\% on energy scheduling, respectively. The gains are larger early in the data stream, a pattern consistent with a benefit from accounting for parameter uncertainty when data are limited.

\end{abstract}

\input{Sections/0_Introduction}
\input{Sections/1_Review}
\input{Sections/2_Framework}
\input{Sections/3_Algorithms}
\input{Sections/4_Experiments}
\input{Sections/5_Conclusions}

\clearpage
\bibliography{aaai2027}

\end{document}

%% file: Sections/0_Introduction.tex
\section{Introduction}
\subsection{Problem Statement}

This paper studies online \emph{contextual optimization under uncertainty}~\citep{Sadana-2025-EJOR}. At each stage $t \in \{1,\ldots,T\}$, the decision-maker observes a context $x_t \in \mathcal{X}$, chooses $z_t \in \mathcal{Z}$, and then observes an uncertain outcome $y_t \in \mathcal{Y}$. Let $\mathcal{H}_{t-1}=\{(x_\tau,y_\tau)\}_{\tau=1}^{t-1}$ denote the history available before stage $t$, and let $\pi_t:\mathcal{X}\rightarrow\mathcal{Z}$ the stage-$t$ policy learned from $\mathcal{H}_{t-1}$, so $z_t=\pi_t(x_t)$. The objective is to minimize the cumulative cost:
\begin{equation}
    \min_{\pi_1,\ldots,\pi_T} \quad \sum_{t = 1}^{T} c\,(\pi_t(x_t); y_t), \label{problem}
\end{equation}
where $c:\mathcal{Z}\times\mathcal{Y}\rightarrow\mathbb{R}$ is a known cost function. We assume full-information feedback: after observing $y_t$, the decision-maker can evaluate $c(z;y_t)$ for every $z\in\mathcal{Z}$~\citep{onlineDFL}. Under bandit feedback, by contrast, only the cost of the chosen action is observed~\citep{Chen-2021-JASA}. This setting captures repeated decisions in inventory management~\citep{Acimovic-2015-MSOM, Wang-2017-POM, Cheung-2022-MS}, transportation~\citep{Jia-2025-TS}, revenue management~\citep{Ferreira-2016-MSOM}, and advertising~\citep{Ye-2023-MS}. Because data arrive sequentially and may evolve over time, policies must adapt as new observations become available~\citep{Hoi-2021-NeuroCom}.

A common strategy to solve Problem~\eqref{problem} combines predict-then-optimize (PtO) with online learning. In PtO, a predictive model estimates $y_t$ from $x_t$, and an optimization model maps that estimate to a feasible decision~\citep{Bertsimas-2020-MS}. PtO accommodates complex feasible regions, while online updates adapt the model to streaming data~\citep{Shalev-2025-Book}. However, minimizing prediction error need not minimize downstream cost because the two objectives may be misaligned~\citep{Liu-2022-arxiv, Kong-2022-NIPS}, particularly under model misspecification~\citep{SPO,SPO2}. Decision-focused learning (DFL) addresses this mismatch by training predictive models directly on downstream decision loss~\citep{DFLreview}. Its theoretical properties and empirical performance have been studied extensively in batch settings~\cite{Balghiti-2023, Geng-2024-NIPS}.

Recent work has extended DFL to online settings~\citep{Liu-2022-arxiv, onlineDFL}, but two limitations remain. First, existing approaches are predominantly frequentist and maintain a single parameter estimate. Ignoring parameter uncertainty, especially early in the data stream, can produce overconfident, high-variance decisions~\citep{Draper-1995-JRSSB, Depeweg-2018-ICML}. Second, most DFL methods require gradients through the optimization problem or problem-specific differentiable surrogate losses, limiting their applicability to generic nondifferentiable problems.

To address both limitations, we propose Bayesian online contextual optimization (\texttt{BOCO}), a decision-focused Bayesian framework that maintains and sequentially updates a belief distribution over model parameters. By propagating this distribution through the downstream optimization problem, \texttt{BOCO} accounts for parameter uncertainty when prescribing decisions. This Bayesian formulation offers two key advantages. First, the generalized posterior can be updated recursively: the current belief summarizes the retained decision evidence and provides the starting distribution for the next update~\citep{Jones-2024-nips, Haddouche-2022-NIPS}. Second, and more importantly, we introduce a decision-focused belief update: we transform historical decision loss into a \emph{generalized power likelihood}~\citep{Holmes-2017-Biom}, so that the posterior assigns greater mass to parameters that induce better decisions. Consequently, the update aligns learning with the downstream objective and does not require specifying a potentially misspecified likelihood for the data-generating process.

Although this formulation is conceptually straightforward, the resulting posterior is generally unavailable in closed form because the likelihood depends on optimization-induced decision loss. We therefore develop two particle-based inference schemes. The first is a sequential Monte Carlo (SMC) sampler~\citep{Dai-2022-JASA}, which represents the posterior with weighted particles and updates them through sequential reweighting and resample-move steps. Its gradient-free Markov chain Monte Carlo (MCMC) rejuvenation makes it applicable to nondifferentiable optimization problems. The second is a function-space Stein variational gradient descent (FSVGD) algorithm~\citep{SVGD, FSVGD}, which represents the posterior with equally weighted particles and evolves them by deterministic gradient flow. The two schemes are complementary: SMC does not require posterior gradients, whereas FSVGD uses them when they are available.

To evaluate both computational regimes, we conduct numerical experiments on two online contextual optimization tasks with real-world data: a discrete knapsack task handled without gradients and an energy-scheduling task for which (sub)gradients through the optimization layer are available. We compare \texttt{BOCO} with frequentist online DFL baselines~\citep{onlineDFL}, reporting the temporal mean and standard deviation of regret within each run, together with percentile bands across independent runs. Across both tasks, BOCO reduces mean regret and its temporal variability relative to the baselines.

In summary, this work makes the following contributions:
\begin{enumerate}
    \item To the best of our knowledge, we introduce BOCO, the first decision-focused Bayesian framework for online contextual optimization. BOCO uses a decision-focused power likelihood to update parameter beliefs and prescribes decisions under the resulting posterior predictive distribution, thereby aligning online learning with downstream decision-making while explicitly accounting for parameter uncertainty.

    \item We develop two particle-based algorithms for tracking the evolving decision-focused posterior. The SMC sampler supports generic nondifferentiable optimization problems through gradient-free MCMC rejuvenation without requiring handcrafted differentiable surrogate losses. FSVGD uses posterior-gradient information to update particles in high-dimensional Bayesian neural network (BNN) models.

    \item We evaluate BOCO on real-world knapsack and energy-scheduling tasks. At the full horizon, SMC reduces mean regret by 6.9\% on knapsack, while FSVGD reduces it by 49.3\% on energy scheduling, relative to the best frequentist baseline. The reductions in both mean regret and temporal regret variability are larger at the early horizon in both tasks, a pattern consistent with a benefit from accounting for parameter uncertainty when data are limited.
\end{enumerate}

%% file: Sections/1_Review.tex
\section{Related Work}
This work connects three strands of research: offline PtO and DFL, their online extensions, and Bayesian methods for sequential decision-making.

\subsubsection{Offline PtO and DFL}
PtO first learns a predictive model from contextual data and then uses its predictions as inputs to a downstream optimization problem. Training the model with a conventional statistical loss, such as mean squared error or negative log-likelihood, can nevertheless yield poor decisions because predictive accuracy need not align with downstream decision quality~\citep{Donti-2017-NIPS, Kong-2022-NIPS}. DFL addresses this target mismatch by training the predictive model end-to-end on the downstream decision loss~\citep{DFLreview}. A central computational challenge is differentiating through the optimization problem; existing work has therefore developed differentiable optimization layers and efficient implicit-gradient methods that integrate with gradient-based training~\citep{Bradon-2017-ICML, Pan-2024-NIPS}.

These approaches primarily produce a single frequentist estimate and typically require either a differentiable optimization problem or a problem-specific surrogate gradient. In contrast, \texttt{BOCO} maintains a belief distribution over model parameters and updates it through a generalized posterior defined by the historical decision loss associated with each parameter. This formulation both represents parameter uncertainty and separates the definition of the posterior from the availability of gradients. Once the posterior is specified, established approximate Bayesian inference methods can be used to track it~\citep{Martin-2024-SS}; in particular, gradient-free MCMC methods allow \texttt{BOCO} to handle nondifferentiable problems without constructing handcrafted differentiable surrogates.

\subsubsection{Online PtO and DFL}
Online PtO extends the offline pipeline by updating the predictive model as observations arrive. In the full-information setting, these updates commonly use online gradient descent (OGD)~\citep{Zinkevich-2003-ICML} or mirror-descent variants. For example, \citet{Lobos-2021-ICML} developed a dual mirror-descent framework for joint online learning and optimization, although its learning objective is not explicitly aligned with downstream decision quality. For a contextual knapsack problem with uncertain rewards and resource consumption, \citet{Liu-2022-arxiv} incorporated the \texttt{SPO+} loss~\citep{SPO} into dual mirror descent to obtain decision-aware model updates.

The closest work to ours is online DFL~\citep{onlineDFL}, which develops decision-focused OGD and follow-the-perturbed-leader (FTPL)~\citep{Suggala-2020-ICALT} algorithms for online contextual linear programs. Although these methods provide theoretical guarantees, they represent the predictive model with a single parameter estimate and therefore do not explicitly propagate parameter uncertainty into the decision. \texttt{BOCO} instead maintains a sequentially updated parameter distribution and prescribes decisions under the induced posterior predictive distribution. This distinction is particularly relevant early in the data stream, when parameter uncertainty is substantial and point estimates may lead to unstable decisions.

\subsubsection{Bayesian Sequential Decision-Making}
Bayesian methods have a long history in sequential decision-making, including bandit problems~\citep{Agrawal-2012-ACLT} and reinforcement learning~\citep{Ghavamzadeh-2015-Book}. Whereas those settings emphasize exploration under partial feedback, we study full-information contextual optimization and use Bayesian inference to represent uncertainty about the predictive model. Structurally, the \texttt{BOCO} policy resembles Bayesian model averaging: it aggregates the predictive scenarios induced by a posterior over model parameters and optimizes against the resulting posterior predictive distribution.

More closely related, \citet{shapiro} developed a Bayesian formulation for stochastic control with side information, using posterior predictive uncertainty to construct decisions. The key distinction lies in how the posterior is updated. Conventional Bayesian approaches update the posterior using a specified probabilistic likelihood, which may be sensitive to misspecification~\citep{Nott-2023-ARSIA, Grunwald-2017-BA}. By contrast, \texttt{BOCO} uses a decision-focused power likelihood derived from historical decision loss. Its belief update therefore rewards parameters for prescriptive performance rather than predictive goodness-of-fit, directly aligning Bayesian learning with the downstream optimization objective.

%% file: Sections/2_Framework.tex
\section{Bayesian Online Contextual Optimization}
We now formalize the online contextual optimization problem and introduce \texttt{BOCO}. Let $T$ be the number of stages, let $\R^d$ denote the $d$-dimensional real vector space, and let $d_x$, $d_y$, $d_z$, and $d_\theta$ denote the dimensions of the context, outcome, decision, and model-parameter vectors, respectively. At each stage $t\in\{1,\ldots,T\}$, the decision-maker observes a context $x_t\in\mathcal{X}\subseteq\R^{d_x}$, selects a decision $z_t\in\mathcal{Z}\subseteq\R^{d_z}$, and then observes an outcome $y_t\in\mathcal{Y}\subseteq\R^{d_y}$. The decision incurs the known cost $c(z_t;y_t)$, where $c:\mathcal{Z}\times\mathcal{Y}\rightarrow\R$. We use $X_t$ and $Y_t$ for the random variables whose realizations are $x_t$ and $y_t$, respectively.

Let $\mathcal{H}_{t-1}\coloneqq\{(x_\tau,y_\tau)\}_{\tau=1}^{t-1}$ denote the observations available before stage $t$. A stage-wise policy $\pi_t:\mathcal{X}\rightarrow\mathcal{Z}$ maps the current context to the decision $z_t=\pi_t(x_t)$ and is learned using $\mathcal{H}_{t-1}$. Thus, the sequence of policies seeks to minimize the cumulative cost in Equation~\eqref{problem}. For notational simplicity, we suppress the dependence of $\pi_t$ on $\mathcal{H}_{t-1}$; the history enters through the learned model parameters, and $x_t$ may also contain history-derived features.

Let $\Theta\subseteq\R^{d_\theta}$ be the parameter space and
\[
    \mathcal{F}_\Theta
    \coloneqq
    \{f(\cdot;\theta):\mathcal{X}\rightarrow\mathcal{Y}
    \mid\theta\in\Theta\}
\]
the predictive-model family. Here, $f(x;\theta)$ is the outcome predicted for context $x$ by parameter vector $\theta$. Before observing any data, we assign $\theta$ an initial probability density $\rho_1^-(\theta)$ on $\Theta$. We assume a continuous parameter space for exposition; integrals can be replaced by sums for a discrete parameter space. At stage $t$, $\rho_t^-$ denotes the parameter density used to make the decision after $x_t$ is observed, whereas $\rho_t^+$ denotes the density after $(x_t,y_t)$ is incorporated. Sequential updating gives $\rho_{t+1}^-=\rho_t^+$.
% ====================================================================
\subsection{Decision-Focused Power Likelihood}
For comparison, suppose that each parameter $\theta$ also specifies a conditional probability model for the outcome. We write
$p_\theta(y\mid x)$ for the conditional probability density of model $f(\cdot;\theta)$ assigned to outcome $y$ given context $x$. Once $(x_t,y_t)$ is observed, $p_\theta(y_t\mid x_t)$ is viewed as a function of $\theta$ and is the stage-$t$ likelihood contribution. The contexts are treated as given, so their distribution is not modeled.

Assuming that the outcomes are conditionally independent across stages given their contexts and $\theta$, the likelihood of
$\mathcal{H}_t=\{(x_\tau,y_\tau)\}_{\tau=1}^{t}$ is
\begin{equation}
    L_p(\mathcal{H}_t\mid\theta)
    \coloneqq
    \prod_{\tau=1}^{t}p_\theta(y_\tau\mid x_\tau).
\end{equation}
Here, $L_p(\mathcal{H}_t\mid\theta)$ denotes the joint likelihood of all outcomes observed through stage $t$, evaluated at parameter $\theta$. Bayes' rule can then be written in either batch or recursive form:
\begin{align}
    \rho_t^+(\theta)
    &=
    \frac{L_p(\mathcal{H}_t\mid\theta)\rho_1^-(\theta)}
    {\int_\Theta L_p(\mathcal{H}_t\mid\nu)\rho_1^-(\nu)\,\mathrm{d}\nu}
    \notag\\
    &=
    \frac{p_\theta(y_t\mid x_t)\rho_t^-(\theta)}
    {\int_\Theta p_\nu(y_t\mid x_t)\rho_t^-(\nu)\,\mathrm{d}\nu}.
\end{align}
Conventional Bayesian updating therefore evaluates parameters according to predictive goodness-of-fit. However, if the family $p_\theta(y\mid x)$ cannot represent the data-generating process, likelihood misspecification may lead to misleading posterior inference~\citep{Nott-2023-ARSIA}. Power likelihoods provide a mechanism for calibrating the influence of the data by replacing $L_p(\mathcal{H}_t\mid\theta)$ with $L_p(\mathcal{H}_t\mid\theta)^\lambda$ for some constant $\lambda>0$~\citep{Ibrahim-2000-SS, Bhattacharya-2019-AOS}.

\texttt{BOCO} instead uses generalized Bayesian updating~\citep{Bissiri-2016-JRSSB} to evaluate parameters by the quality of their induced decisions. For any predicted outcome $\hat y\in\mathcal{Y}$, define
\[
    z^*(\hat y)
    \in
    \argmin_{z\in\mathcal{Z}}c(z;\hat y),
\]
where $\hat y$ is treated as the outcome when solving the downstream optimization problem. If the minimizer is not unique, $z^*(\hat y)$ denotes any fixed selection from the set of minimizers. The stage-wise decision loss and its cumulative value are
\begin{align}
    \ell_t(\theta)
    &\coloneqq
    c\lra{z^*\lra{f\lra{x_t;\theta}};y_t},
    \label{stage_loss}\\
    V_t(\theta)
    &\coloneqq
    \sum_{\tau=1}^{t}\ell_\tau(\theta).
    \label{historical_loss}
\end{align}
Thus, $\ell_t(\theta)$ evaluates the decision that parameter $\theta$ would induce from $x_t$ against the subsequently observed outcome $y_t$. Accordingly, $V_t(\theta)$ is the cumulative counterfactual decision loss associated with using that parameter throughout the observed history.

\texttt{BOCO} converts this additive loss into the multiplicative generalized likelihood:
\begin{equation}
    L(\mathcal{H}_t\mid\theta)
    \coloneqq
    \exp\{-V_t(\theta)\}
    =
    \prod_{\tau=1}^{t}\exp\{-\ell_\tau(\theta)\}.
    \label{df-likelihood}
\end{equation}
The quantity $L$ is not required to be a probability density over the data; rather, it is a nonnegative weight that favors parameters with smaller historical decision loss. In particular, it recovers the conventional likelihood when the stage loss coincides with the negative log-likelihood~\citep{Germain-2016-NIPS}:
\begin{equation*}
    \ell_t(\theta) = -\log p_\theta(y_t\mid x_t) \Longrightarrow L(\mathcal{H}_t\mid\theta)
    =
    \prod_{\tau=1}^{t}p_\theta(y_\tau\mid x_\tau).
\end{equation*}
This special case illustrates that conventional Bayesian updating is contained within the loss-based construction.

Because the numerical scale of the loss determines how strongly the data concentrate the generalized posterior, an untempered likelihood may produce an overly concentrated or overly diffuse belief~\citep{Holmes-2017-Biom, Wu-2023-BA}. We therefore introduce the decision-focused power likelihood:
\begin{equation}
    L_\lambda(\mathcal{H}_t\mid\theta)
    \coloneqq
    \exp\{-\lambda_tV_t(\theta)\},
\end{equation}
where $\lambda_t>0$ is the stage-$t$ learning rate, also called the inverse temperature. A larger $\lambda_t$ gives the accumulated decision loss greater influence relative to the initial prior, whereas a smaller value gives the prior greater influence. Unlike a fixed power $\lambda$, $\lambda_t$ is allowed to vary with the amount of data observed.

\subsection{Decision-Focused Online Belief Update}
Combining the power likelihood with the initial prior gives the decision-focused posterior:
\begin{equation}\label{dfp-update}
    \rho_t^+(\theta) \propto \exp\{-\lambda_tV_t(\theta)\}\rho_1^-(\theta)
\end{equation}
Parameters with smaller $V_t(\theta)$ receive greater posterior density, subject to the initial prior $\rho_1^-$. Unlike frequentist DFL, which represents the learned model through a point estimate, $\rho_t^+$ retains a distribution over parameter values and thereby represents parameter uncertainty.

Although Equation~\eqref{dfp-update} uses the full history, it admits a recursive form. Because $V_t(\theta)=V_{t-1}(\theta)+\ell_t(\theta)$, for $t\geq2$,
\begin{equation}\label{sequential}
    e^{-\lambda_tV_t(\theta)} = e^{-\lambda_t\ell_t(\theta) -(\lambda_t-\lambda_{t-1})V_{t-1}(\theta)} e^{-\lambda_{t-1}V_{t-1}(\theta)}.
\end{equation}
Substitution into Equation~\eqref{dfp-update} yields
\begin{equation}
    \rho_t^+(\theta)
    \propto
    e^{-\lambda_t\ell_t(\theta) -(\lambda_t-\lambda_{t-1})V_{t-1}(\theta)}
    \rho_{t-1}^+(\theta).
    \label{sequentialweight}
\end{equation}
The correction $(\lambda_t-\lambda_{t-1})V_{t-1}(\theta)$ is required because the inverse temperature may change with $t$; if $\lambda_t=\lambda_{t-1}$, the update reduces to multiplication by $e^{-\lambda_t\ell_t(\theta)}$. 

The learning rate controls the tradeoff between the initial prior and the accumulated decision evidence. A large $\lambda_t$ concentrates the posterior near parameters that minimize $V_t(\theta)$ and may underestimate parameter uncertainty. Conversely, a small value retains more prior influence and slows posterior concentration. Motivated by the inverse-square-root schedule used in PAC-Bayes\footnote{PAC stands for \emph{probably approximately correct}.} analysis~\citep{Alquier-2024-FTML}, for a stationary stream we set $\lambda_t=\lambda_1t^{-1/2}$
where $\lambda_1>0$ is the user-specified initial learning rate.

For a nonstationary stream, older observations may no longer represent the current environment. We therefore introduce a window width $H\in\mathbb{N}_+$ and define
\begin{align*}
    \tau_t \coloneqq \max\{1,t-H+1\}, \quad \lambda_t = \lambda_1\bigl(\min\{t,H\}\bigr)^{-1/2}, \notag\\
    V_t^{H}(\theta) \coloneqq \sum_{\tau=\tau_t}^{t}\ell_\tau(\theta),\quad  L_\lambda^{H}(\mathcal{H}_t\mid\theta) \coloneqq e^{-\lambda_tV_t^{H}(\theta)}.
\end{align*}
Here, $\tau_t$ is the first stage retained in the current window,
$V_t^{H}(\theta)$ is the decision loss over the most recent
$\min\{t,H\}$ observations, and $L_\lambda^{H}$ is the associated windowed power likelihood. The schedule for $\lambda_t$ decreases until the window contains $H$ observations and then remains at
$\lambda_1H^{-1/2}$. When the sliding window is used, $V_t(\theta)$ in Equation~\eqref{dfp-update} is replaced by $V_t^{H}(\theta)$. The recursion in Equations~\eqref{sequential}--\eqref{sequentialweight} applies to the untruncated cumulative loss. Under the sliding window, the implementation instead updates $V_t^{H}$ by adding $\ell_t$ and, once $t>H$, removing the exiting loss $\ell_{t-H}$.

This sliding-window scheme differs from geometric power weighting, which assigns observation $\tau$ the stage-dependent weight
$\lambda_{\tau,t}=\lambda_1\alpha^{t-\tau}$ for a retention factor $\alpha\in(0,1)$~\citep{McCarthy-2016-AOAS}. The geometric scheme gradually discounts older observations. In contrast, our window assigns the common weight $\lambda_t$ to observations with
$\tau\in\{\tau_t,\ldots,t\}$ and zero weight to earlier observations. In addition, the sliding window reduces the number of downstream optimization problems required to evaluate a newly generated particle from $t$ to at most $H$, as discussed in Section~\ref{sec:4}.

\subsection{Bayesian Ensemble Policy}
At stage $t$, the density $\rho_t^-$ summarizes observations through stage $t-1$ and represents uncertainty about the model parameter before $y_t$ is observed. For a context $x\in\mathcal{X}$ and a measurable set
$A\subseteq\mathcal{Y}$, define the posterior predictive distribution
\begin{equation}
    \hat P_t(A\mid x)
    \coloneqq
    \int_\Theta
    \delta_{f(x;\theta)}(A)\rho_t^-(\theta)\,\mathrm{d}\theta.
\end{equation}
Here, $\delta_a$ is the Dirac point-mass distribution at $a$: it assigns probability one to any set containing $a$ and zero otherwise. Equivalently, $\hat P_t(\cdot\mid x)$ is the distribution obtained by drawing
$\theta$ from $\rho_t^-$ and setting the predicted outcome to
$f(x;\theta)$. Thus, it propagates parameter uncertainty into a distribution over predicted outcomes.

The Bayesian ensemble policy minimizes expected cost under this posterior predictive distribution:
\begin{equation}
    \pi_t(x_t)
    \in
    \argmin_{z\in\mathcal{Z}}
    \E_{Y\sim\hat P_t(\cdot\mid x_t)}
    \lrb{c(z;Y)}.
    \label{Policy}
\end{equation}
The notation $Y\sim\hat P_t(\cdot\mid x_t)$ means that the random outcome $Y$ follows the posterior predictive distribution. The expectation in Equation~\eqref{Policy} is equivalently
$\int_{\mathcal{Y}}c(z;y)\hat P_t(\mathrm{d}y\mid x_t)$.

Importantly, Equation~\eqref{Policy} does not generally reduce to optimizing the posterior-mean prediction
\begin{equation*}
    \bar y_t
    \coloneqq
    \int_\Theta
    f(x_t;\theta)\rho_t^-(\theta)\,\mathrm{d}\theta,
\end{equation*}
where $\bar y_t$ is the average prediction under the current parameter belief. The two formulations are guaranteed to be equivalent when \(c(z;y)\) is affine in \(y\) for every fixed \(z\), because then
\[
    \E[c(z;Y)]
    =
    c\lra{z;\E[Y]}.
\]
For a nonlinear cost, replacing the posterior predictive distribution with its mean discards distributional information and can lead to uncertainty-agnostic decisions. In our experiments, the knapsack cost is affine in $y$, whereas the energy-scheduling cost is nonlinear in \(y\), so replacing the predictive distribution with its mean is not generally valid; we therefore retain the empirical predictive distribution.

%% file: Sections/3_Algorithms.tex
\section{Inference Algorithms}\label{sec:4}
The posterior in Equation~\eqref{dfp-update} is generally unavailable in closed form because $V_t(\theta)$ depends on the predictive model and the downstream optimizer. We therefore develop two nonparametric particle approximations to the posterior: SMC, which requires no gradients, and FSVGD, which uses posterior-gradient information when it is available. Both algorithms use the same optimization stage but different learning updates. The complete procedures of the two algorithms are provided in the supplementary material.

\subsection{Optimization Stage}
At stage $t$, both algorithms maintain
$\mathcal{P}_t=\lrc{(w_t^i,\theta_t^i)}_{i=1}^{n}$ to approximate $\rho_t^-$, where $w_t^i\geq0$ and $\sum_iw_t^i=1$; FSVGD uses $w_t^i=1/n$. The resulting empirical posterior predictive distribution and policy are
\begin{align*}
    \hat P_t(\cdot\mid x_t)
    &=
    \sum_{i=1}^{n}w_t^i\delta_{f(x_t;\theta_t^i)}(\cdot),\\
    \pi_t(x_t)
    &\in
    \argmin_{z\in\mathcal{Z}}
    \sum_{i=1}^{n}w_t^i c\lra{z;f(x_t;\theta_t^i)}.
\end{align*}
We use $n=32$ in all experiments; increasing $n$ trades additional computation for a more accurate particle approximation.

\subsection{Learning in SMC}
After observing $y_t$, SMC first reweights the particles using the sliding-window analogue of Equation~\eqref{sequentialweight}. For $t\geq2$,
\begin{align*}
    \tilde w_{t+1}^i
    &=
    w_t^i
    e^{-\lambda_t(\ell_t(\theta_t^i) - \ell_{t-H}(\theta_t^i)) -(\lambda_t-\lambda_{t-1})V^H_{t-1}(\theta_t^i) },\\
    w_{t+1}^i
    &=
    \frac{\tilde w_{t+1}^i}
    {\sum_{j=1}^{n}\tilde w_{t+1}^j}.
\end{align*}
At $t=1$, the update follows $\tilde{w}_2^i=w^i_1e^{-\lambda_1\ell_1(\theta^i_1)}$. In addition, $\ell_{t-H} = 0$ if $t \leq H$.

Repeated reweighting can concentrate most mass on a few particles. Following the resample-move strategy~\citep{Arulampalam-2002-IEEESP, Gilks-2001-JRSSB}, we monitor the effective sample size
\begin{equation}
    n_{\mathrm{eff}}
    \coloneqq
    \frac{1}{\sum_{i=1}^{n}(w_{t+1}^i)^2}.
\end{equation}
If $n_{\mathrm{eff}}\leq 0.5n$, we sample ancestor indices:
\begin{equation*}
    a_i
    \sim
    \operatorname{Categorical}
    (w_{t+1}^1,\ldots,w_{t+1}^n),
    \qquad
    \bar\theta_t^i
    =
    \theta_t^{a_i},
\end{equation*}
and reset the weights $w_{t+1}^i = 1/n$.  

Resampling removes negligible-weight particles but may duplicate high-weight ancestors. 
We restore diversity with elliptical slice sampling~\citep{Murray-2010-AISTATS}, a gradient-free and tuning-free MCMC method suited to the Gaussian prior we adopt for the model. Let $\mathcal{K}_t^{\mathrm{ESS}}$ be its transition kernel, which leaves $\rho_t^+$ invariant:
\begin{equation*}
    \int_\Theta
    \rho_t^+(\mathrm{d}\theta)
    \mathcal{K}_t^{\mathrm{ESS}}(\theta,\mathrm{d}\theta')
    =
    \rho_t^+(\mathrm{d}\theta').
\end{equation*}
Then, each resampled particle is rejuvenated as
\begin{equation*}
    \theta_{t+1}^i
    \sim
    \mathcal{K}_t^{\mathrm{ESS}}
    (\bar\theta_t^i,\cdot).
\end{equation*}
In implementation, we conduct $K>1$ steps of ESS moves for each particle for stationarity. If resampling is not triggered, the particles retain their locations and updated weights. Evaluating an elliptical-slice transition requires the windowed power likelihood $L^H_\lambda$ and hence the windowed cumulative loss $V_t^H$, which is cached for subsequent reweighting.

\subsection{Learning in FSVGD}
FSVGD transports equally weighted particles using the posterior score
$\nabla_\theta\log\rho_t^+(\theta)$ and does not impose a parametric variational family. To approximate a probability distribution $\rho$ using a particle set $\{\theta_1,\dots,\theta_n\}$, the standard SVGD~\cite{SVGD} uses:
\begin{align*}
    \theta_i & \gets \theta_i+\eta\phi(\theta_i),\\
    \phi(\theta_i) &= \frac{1}{n} \sum_{j=1}^{n} \lra{
        k(\theta_j,\theta_i) \nabla_{\theta_j}\log\rho(\theta_j) + \nabla_{\theta_j}k(\theta_j,\theta_i)
    },
\end{align*}
where $\eta>0$ is the step size and $k$ is a positive-definite kernel. The two terms on the right-hand side are the posterior attraction and the repulsion, respectively.

Based on the standard SVGD, we propose the online function-space update:
\begin{equation}\label{onlineFSVGD}
    \begin{aligned}
        \theta_{t+1}^i
        &\gets
        \theta_t^i+\eta\phi_t^i,\\
        \phi_t^i
        &=
        \frac{1}{n}
        \sum_{j=1}^{n}
        \lra{
            k_\mathcal{Y}(\theta_t^i,\theta_t^j)
            \nabla_{\theta_t^j}
            \log\rho_t^+(\theta_t^j)
            -
            \nabla_{\theta_t^i}
            k_\mathcal{Y}(\theta_t^i,\theta_t^j)
        }.
    \end{aligned}
\end{equation}
Relative to standard SVGD, Equation~\eqref{onlineFSVGD} uses a function-space kernel and a computationally cheaper repulsion.

Following FSVGD~\citep{FSVGD}, we measure particle similarity through their predictions on the observed contexts:
\begin{align*}
    d_t(\theta_t^i,\theta_t^j)
    &\coloneqq
    \frac{1}{\min\{t,H\}}
    \sum_{\tau=\tau_t}^{t}
    \left\|
        f(x_\tau;\theta_t^i)
        -
        f(x_\tau;\theta_t^j)
    \right\|_2^2,\\
    k_\mathcal{Y}(\theta_t^i,\theta_t^j)
    &\coloneqq
    \exp\lrc{
        -d_t(\theta_t^i,\theta_t^j)/h_t
    }.
\end{align*}
This kernel can distinguish functionally different predictors even when parameter-space distance is uninformative. We set its bandwidth using the median heuristic~\citep{SVGD}:
\begin{equation*}
    h_t
    =
    \frac{
        \operatorname{median}_{i<j}
        d_t(\theta_t^i,\theta_t^j)
    }{\log(n+1)}.
\end{equation*}

In addition, standard SVGD differentiates the repulsion through each $\theta_t^j$. We instead use $-\nabla_{\theta_t^i}\sum_jk_\mathcal{Y}(\theta_t^i,\theta_t^j)/n$, treating the other particles as fixed.
In our implementation, this approximation reduces the number of model-Jacobian evaluations required for the repulsive term. This substitution is a computational heuristic and does not inherit the theoretical guarantees of standard SVGD.

Finally, the required posterior score is
\begin{equation*}
    \nabla_{\theta_t^i} \log\rho_t^+(\theta_t^i) = \nabla_{\theta_t^i} \lra{ -\lambda_tV^H_t(\theta_t^i) + \log\rho_1^-(\theta_t^i) }.
\end{equation*}
Because the sliding window retains at most $H$ recent observations, computing $V^H_t$ and its gradient requires at most $H$ differentiations of the optimization layer per particle. In implementation, we conduct $K>1$ gradient steps for each particle.

%% file: Sections/4_Experiments.tex
\section{Numerical Experiments}
We evaluate \texttt{BOCO} on two online contextual optimization tasks with real-world data: a discrete knapsack task handled without gradients and an energy-scheduling task for which (sub)gradients through the optimization layer are available. The two tasks evaluate the framework in complementary computational regimes and test whether posterior averaging improves decision quality and stability when data are limited. We compare against the decision-focused OGD and FTPL methods of~\citet{onlineDFL}, using SMC for knapsack and FSVGD for energy scheduling. We run all experiments in Python 3.12 and PyTorch 2.11 with two Intel Xeon Gold 6230 processors. Model architectures, hyperparameters, and other implementation details are provided in the supplementary material. 

We conduct 10 independent runs of each approach with random initializations and compute the running mean regret \(\bar R_t\) and temporal standard deviation of regret \(\sigma_t\) within each run as follows:
\begin{align*}
    r_t       &= c(\pi_t(x_t);y_t) - \min_{z\in\mathcal{Z}} c(z;y_t), \\
    \bar{R}_t &= \frac{1}{t-1} \sum_{\tau=2}^{t} r_\tau, \quad \sigma_t = \sqrt{\frac{1}{t-2}\sum_{\tau=2}^{t}(r_\tau - \bar{R}_t)^2}.
\end{align*}
We compute both metrics starting from $\tau = 2$ because all policies are randomly initialized at $\tau=1$. For each metric, Tables~\ref{tab:1} and~\ref{tab:2} further report the mean and standard deviation across 10 independent runs.

\subsection{Knapsack}
Consider a knapsack problem with uncertain item values
\begin{equation*}
    \min_{z \in \{0,1\}^{48}} \, c(z;y_t) = -y_t^\top z, \quad \text{s.t.}\quad W^\top z \leq 60,
\end{equation*}
where $W\in\R^{48}$ is the fixed item-weight vector provided in the supplementary material. Because $z$ is binary, the solution mapping $z^*(y)$ is nondifferentiable with respect to the value vector $y$. Consequently, OGD and FTPL use the log-barrier surrogate proposed by~\citet{onlineDFL}, whereas SMC operates without differentiating through the knapsack problem.

We follow the knapsack-electricity (KE) experiment of~\citet{Geng-2024-NIPS}, which treats the 48 half-hourly electricity prices within a day as item values. We adapt the experiment to online learning in three ways. First, we normalize each price dimension by its mean to remove the fixed intraday scale pattern. Second, we use a neural network with two hidden layers of 32 units, yielding $d_\theta=1{,}377$ parameters. Third, rather than constructing offline training and test sets, we retain the observations in their original temporal order.

\begin{figure}[t]
    \centering
    \includegraphics[width=\linewidth]{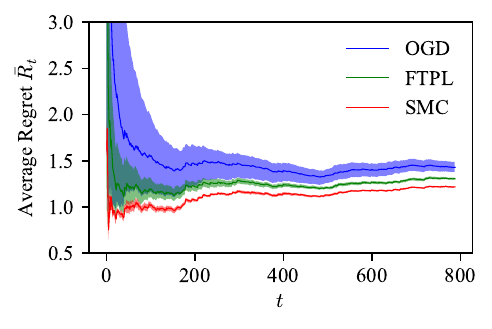}
    \caption{Average regret $\bar R_t$ for the knapsack problem. Lines show the mean across 10 independent runs, and shaded bands show the 10th--90th percentiles across runs at each stage.}
    \label{fig:1}
\end{figure}

Figure~\ref{fig:1} shows that SMC achieves lower $\bar{R}_t$ and a narrower cross-run percentile band than both frequentist baselines. OGD exhibits the greatest dispersion, particularly near the beginning of the stream, whereas SMC stabilizes more quickly. This behavior is consistent with the proposed benefit of accounting for parameter uncertainty rather than relying on a single parameter estimate when only limited data are available.

Table~\ref{tab:1} reports $\bar R_t$ and $\sigma_t$, averaged across the 10 runs, at the early horizon $t=100$ and the full horizon $t=789$. 
Relative to FTPL, SMC reduces $\bar R_{100}$ by $13.6\%$ and $\sigma_{100}$ by $24.1\%$. At $t=789$, the corresponding reductions are $6.9\%$ and $6.3\%$. Thus, both reductions are larger at the early horizon, while SMC remains better than FTPL on both metrics at the full horizon.
% \begin{table}[t]
%     \centering
%     \small
%     \setlength{\tabcolsep}{3.5pt}
%     \begin{tabular}{lcccc}
%         \toprule
%         Method
%         & $\bar R_{100}$
%         & $\sigma_{100}$
%         & $\bar R_{789}$
%         & $\sigma_{789}$\\
%         \midrule
%         OGD  & $1.55 \pm0.47$ & $1.48 \pm 0.70$ & $1.43 \pm0.07$ & $1.33 \pm 0.13$ \\
%         FTPL & $1.18 \pm0.09$ & $1.16 \pm 0.34$ & $1.31 \pm0.01$ & $1.26 \pm 0.05$ \\
%         SMC  & $1.03 \pm0.03$ & $0.89 \pm 0.06$ & $1.22 \pm0.01$ & $1.18 \pm 0.01$ \\
%         \bottomrule
%     \end{tabular}
%     \caption{Knapsack performance at the early ($t=100$) and full horizon ($t=789$). Each entry corresponds to the average metric and one standard deviation computed from 10 independent runs. Lower values are better.}
%     \label{tab:1}
% \end{table}

\begin{table}[t]
    \centering
    \small
    \setlength{\tabcolsep}{3.5pt}
    \begin{tabular}{lcccc}
        \toprule
        Method
        & $\bar R_{100}$
        & $\sigma_{100}$
        & $\bar R_{789}$
        & $\sigma_{789}$\\
        \midrule
        OGD & $1.55 \pm0.47$ & $1.48 \pm 0.70$ & $1.43 \pm0.07$ & $1.33 \pm 0.13$ \\
        FTPL & $1.18 \pm0.09$ & $1.16 \pm 0.34$ & $1.31 \pm0.01$ & $1.26 \pm 0.05$ \\
        \rowcolor{gray!20}
        SMC & $1.02 \pm0.03$ & $0.88 \pm 0.05$ & $1.22 \pm0.01$ & $1.18 \pm 0.01$ \\
        \bottomrule
    \end{tabular}
    \caption{Knapsack performance at the early ($t=100$) and full horizon ($t=789$). Each entry reports the mean \(\pm\) standard deviation of the per-run metric across 10 independent runs. Lower values are better; the shaded row contains the best results.}
    \label{tab:1}
\end{table}
\subsection{Energy Scheduling}
We next consider the day-ahead energy-scheduling problem of~\citet{Donti-2017-NIPS}. The decision $z\in\R_+^{24}$ specifies hourly energy generation, and $y\in\R_+^{24}$ is the realized demand:
\begin{align*}
    \min_{z\geq0}
    \quad&
    c(z;y)
    =
    \frac{1}{24}
    \sum_{s=1}^{24}
    \lra{
        c_u[y_s-z_s]^+
        +
        c_o[z_s-y_s]^+
    },\\
    \text{s.t.}
    \quad&
    |z_{s+1}-z_s|\leq0.4,
    \qquad s=1,\ldots,23.
\end{align*}
Here, $c_u=8$ and $c_o=2$ are the unit underage and overage costs, respectively, $[a]^+=\max\{0,a\}$, and the constraints limit changes in generation between consecutive hours.

Given particle predictions $y_t^i=f(x_t;\theta_t^i)$ and weights $w_t^i$, the sample average approximation used in our policy is
\begin{align*}
    \min_{z\geq0}
    \quad&
    \sum_{i=1}^{n}
    \frac{w_t^i}{24}
    \sum_{s=1}^{24}
    \lra{
        c_u[y_{t,s}^i-z_s]^+
        +
        c_o[z_s-y_{t,s}^i]^+
    },\\
    \text{s.t.}
    \quad&
    |z_{s+1}-z_s|\leq0.4,
    \qquad s=1,\ldots,23.
\end{align*}
We use the dataset of~\citet{Donti-2017-NIPS} and a predictive network mapping 149-dimensional features to a 24-dimensional nonnegative demand vector, with $d_\theta=78{,}624$ parameters. We retain the first $1{,}000$ observations in their temporal order for this task.

\begin{figure}[t]
    \centering
    \includegraphics[width=\linewidth]{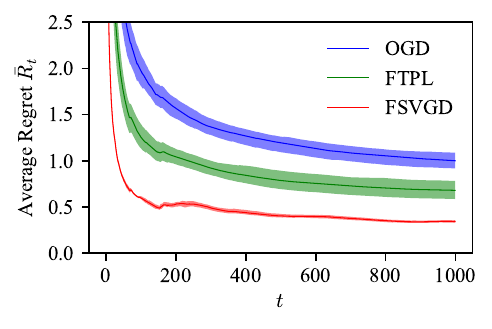}
    \caption{Average regret $\bar R_t$ for the online energy-scheduling problem. Lines show the mean across 10 independent runs, and shaded bands show the 10th--90th percentiles across runs at each stage.}
    \label{fig:2}
\end{figure}

Figure~\ref{fig:2} shows that FSVGD consistently attains lower running mean regret and less cross-run dispersion than OGD and FTPL. Table~\ref{tab:2} further shows that, relative to FTPL, FSVGD reduces \(\bar R_{100}\) by 52.6\% and \(\sigma_{100}\) by 55.2\%. At \(t=1{,}000\), the corresponding reductions are 49.3\% and 38.1\%.

Both the mean-regret and temporal-variability reductions are larger at \(t=100\) than at the full horizon. Nevertheless, FSVGD remains better than FTPL on both metrics at \(t=1{,}000\), indicating that the improvement persists as data accumulate.
This result is particularly relevant to energy scheduling: asymmetric underage and overage costs and ramping constraints couple the 24 hourly decisions, so errors in a single predicted demand trajectory can affect the entire schedule. The posterior-predictive policy evaluates an empirical ensemble of particle-predicted demand trajectories, whereas OGD and FTPL commit to a single-point prediction.

\begin{table}[t]
    \centering
    \small
    \setlength{\tabcolsep}{3.5pt}
    \begin{tabular}{lcccc}
        \toprule
        Method
        & $\bar R_{100}$
        & $\sigma_{100}$
        & $\bar R_{1000}$
        & $\sigma_{1000}$\\
        \midrule
        OGD & $1.86 \pm0.16$ & $1.64 \pm 0.12$ & $0.99 \pm0.07$ & $0.63 \pm 0.05$ \\
        FTPL & $1.14 \pm0.12$ & $1.16 \pm 0.22$ & $0.67 \pm0.11$ & $0.42 \pm 0.07$ \\
        \rowcolor{gray!20}
        FSVGD & $0.54 \pm0.01$ & $0.52 \pm 0.02$ & $0.34 \pm0.01$ & $0.26 \pm 0.01$ \\
        \bottomrule
    \end{tabular}
    \caption{Energy-scheduling performance at the early and full horizons. Each entry reports the mean \(\pm\) standard deviation across 10 independent runs. Lower values are better; the shaded row contains the best results.}
    \label{tab:2}
\end{table}

%% file: Sections/5_Conclusions.tex
\section{Concluding Remarks}
We introduced \texttt{BOCO}, a generalized Bayesian framework for decision-focused online contextual optimization. \texttt{BOCO} updates a distribution over model parameters through a decision-focused power likelihood and prescribes decisions under the induced posterior predictive distribution, thereby aligning online learning with downstream decision quality while accounting for parameter uncertainty. 
For posterior inference, we developed an SMC sampler with gradient-free rejuvenation and an FSVGD algorithm that uses posterior gradients when they are available. Both algorithms use a sliding window intended to facilitate adaptation to nonstationary data.
Across the two tasks, \texttt{BOCO} achieves lower running mean regret, lower temporal regret variability, and narrower cross-run percentile bands than the evaluated baselines. The reduction in running mean regret is larger over the first 100 stages in both tasks, a pattern consistent with a benefit from accounting for parameter uncertainty in small-sample regimes. Future work may establish regret guarantees for the decision-focused online posterior, select the learning rate and window width adaptively, and develop more scalable posterior approximations.